\documentstyle [12pt] {article}

\newcommand {\est}{\exists^{\rm st}}

\newcommand {\fst}{\forall\hspace{0.5pt}^{\rm st}}

\newcommand {\fstf}{\forall\hspace{0.5pt}^{\rm stfin}}

\newcommand {\kpa}{\kappa}

\newcommand {\vpi}{\varphi}
\newcommand {\vt}{\vartheta}

\newcommand {\lra}{\longrightarrow}
\newcommand {\llra}{\longleftrightarrow}

\newcommand {\bV}{{\bf V}}

\newcommand {\bL}{{\bf L}}
\newcommand {\bS}{{\bf S}}
\newcommand {\bN}{\omega}

\newcommand {\cj}{\;\,\&\;\,}

\newcommand {\sV}{\hbox {$\hspace{0mm}^\ast\!V$}}

\newcommand {\sif}{{\rm st}-$\hspace{-1mm}\in\hspace{-1mm}$-formula}

\newcommand {\ins}{\hspace{3pt}^\ast\hspace{-6pt}\in}
\newcommand {\eqs}{\hspace{3pt}^\ast\hspace{-6pt}=}
\newcommand {\sta}{\hspace{0pt}^\ast\!}

\begin{document}
\title{IST is more than an algorithm to prove ZFC theorems}
\author{V.Kanovei\thanks{On leave to Bergische Universit\"at --
Gesamthochschule Wuppertal during Spring 1993.}\vspace{2mm}\\
Moscow Transport Engineering Institute\vspace{-2mm}
\and
Moscow State University}
\date{May 1993}

\maketitle
\begin{abstract}
There is a sentence in the language of ${\sf IST},$ Nelson's
internal set theory, which is not equivalent in ${\sf IST}$
to a sentence in the $\hspace{-1mm}\in\hspace{-1mm}$-language.
Thus the Reduction algorithm, that converts bounded ${\sf IST}$
formulas with standard parameters to provably (in ${\sf IST})$
equivalent $\hspace{-1mm}\in\hspace{-1mm} $-formulas, cannot be
extended to all formulas of the ${\sf IST}$ language.
\end{abstract}

{\bf Introduction.} Internal set theory ${\sf IST}$ was invented
by Nelson [1977] as an attempt to develop nonstandard mathematics
from a unified axiomatical standpoint. This theory has
demonstrated its ability to ground various branches of nonstandard
analysis, see e.g. van den Berg [1987], F. and M.Diener [1988],
F.Diener and Reeb [1989], Reeken [1992].

It is regarded as one of the advantages of ${\sf IST}$ that
there exists a simple algorithm, introduced also by Nelson, to
transform sentences in the language of ${\sf IST}$ to provably
equivalent (in the sense of provability in ${\sf IST})$ sentences
formulated in the ${\sf ZFC}$ language. This algorithm, together
with Nelson's theorem that ${\sf IST}$ is a conservative extension
of ${\sf ZFC},$ is used sometimes (see e.g. Nelson [1988]) to give
back to the statement that ${\sf IST}$ is nothing more than a new
way to investigate the standard ${\sf ZFC}$ universe.

This is true, indeed, so far as {\it bounded\/} ${\sf IST}$ formulas
are considered. (The mentioned algorithm works for these formulas
only.)

It is the aim of this paper to demonstrate that there is a
certain, explicitly given
sentence in the ${\sf IST}$ language which is not provably
equivalent in ${\sf IST}$ to a sentence in the
$\hspace{-1mm}\in\hspace{-1mm}$-language. Thus the ${\sf IST}$ truth
cannot be completely reduced to the ${\sf ZFC}$ truth.

A sentence of this kind has to be undecidable in ${\sf IST};$
actually the sentence we consider belongs to a type of undecidable
sentences discovered and studied in Kanovei [1991]. It is as
follows:\vspace{2mm}

$(\ast) \ \ \forall\,F\;[\,\fst n\,(F(n) {\rm \ is\ standard)\
}\;\longrightarrow\;\est G\;\fst n\;(F(n)=G(n))\,].$\vspace{2mm}

\noindent ($n$ is assumed to range over integers, $F$ and $G$ over
functions defined on integers and taking arbitrary
values.)\vspace{3mm}

{\bf Theorem 1.} {\it Let\/ $\Phi $ be an arbitrary\/
$\hspace{-1mm}\in\hspace{-1mm}$-sentence. Then the equivalence
$\Phi \;\llra\; (\ast)$ is not a theorem of\/ ${\sf IST}$ unless\/
${\sf IST}$ is inconsistent .}\vspace{3mm}

(Take notice that ${\sf ZFC}$ and ${\sf IST}$ are equiconsistent.) The
idea of the proof is to construct a ${\sf ZFC}$ model $V$ which has
two different ${\sf IST}$ extensions, $\sV$ and $\sV ',$ such that
$(\ast)$ is false in $\sV$ but true in $\sV '.$ Both $\sV$ and
$\sV '$ are elementary extensions of $V$ with respect to
$\hspace{-1mm}\in\hspace{-1mm}$-sentences by the ${\sf IST}$
Transfer, hence true (parameterfree)
$\hspace{-1mm}\in\hspace{-1mm}$-sentences are the same in both
extensions. This proves the
theorem. This reasoning is carried out in the assumption of the
existence of a cardinal $\vt $ such that $\bV\!_\vt,$ the $\vt
$th level of the von Neumann hierarchy of sets, is a ${\sf ZFC}$
model. It will be shown at the end of the paper how this
assumption can be abandoned. \vspace{2mm}

{\sf Acknowledgement.} The author is in debt to M.Reeken and
S.Albeverio for their interest to this research direction and
practical help.\vspace{5mm}

{\bf Preliminaries.} Theory ${\sf IST}$ was introduced in Nelson
[1977]. The ${\sf IST}$ language contains, together with equality,
the membership predicate $\in $ and the standardness predicate st.
Formulas of this language are called\break
{\it \sif s} while formulas of the ${\sf ZFC}$ language are called
$\hspace{-0.8mm}\in\hspace{-1mm}$-{\it formulas}, and also {\it
internal\/} formulas. Two abbreviations are very useful: $\est
x\,...$ and $\fst x\,...$ (there exists standard $x$..., for all
standard $x$...).

${\sf IST}$ contains all axioms of ${\sf ZFC}$ (Separation and
Replacement are formulated in the
$\hspace{-0.8mm}\in\hspace{-1mm}$-language) together with the
following three additional principles or (schemes of)
axioms.\vspace{2mm}

\noindent {\sf Idealization I:} $\fstf A\; \exists\,x\;\forall
\,a \in A\;\Phi (x,a)\,\;\llra\,\;\exists\,x\;\fst a\;\Phi
(x,a)$\vspace{-1mm}

\hfill for any internal formula $\Phi (x,a).$\vspace{2mm}

\noindent {\sf Standardization S:} $\fst X\;\est Y\;\fst x\;[\,x
\in Y\,\;\llra\,\;x \in X\cj \Phi (x)\,]$\vspace{-1mm}

\hfill for any \sif\ $\Phi .$\vspace{2mm}

\noindent {\sf Transfer T:} $\exists\,x\;\Phi (x)\,\;\lra\,\;\est
x\;\Phi (x)$\vspace{-1mm}

\hfill for any internal formula $\Phi (x)$ with standard
parameters. $\Box$\vspace{2mm}

\noindent The formula $\Phi $ can, of course, contain {\it
arbitrary\/} parameters in {\sf I} and {\sf S}.\vspace{2mm}

Thus ${\sf IST} ={\sf ZFC} +{\sf I} +{\sf S} + {\sf T}.$ We would
refer to this theory as {\it plain} ${\sf IST}$ since there have
been introduced several modifications, say, by a kind of
superstructure over ${\sf IST}$ (see Nelson [1988]) with the purpose
to extend the Reduction algorithm to a more wide class of \sif s.

Let $V$ be a transitive ${\sf ZFC}$ model in the ${\sf ZFC}$ universe.
We say that an ${\sf IST}$ model $\sV= \langle \sV;\ins,\eqs,\;
\sta{\rm st}\rangle $ is a {\it regular\/} ${\sf IST}$ extension of
$V$ if and only if there exists an 1-1 embedding $^\ast$: $V$
onto a subset of $\sV$ satisfying\vspace{1mm}

1) $x\in y\;\llra\;\sta x\ins \sta y$ \ and \ $x=y\;\llra\;\sta
x\eqs \sta y$ \ \ for all $x,y\in V,$ and

2) $\sta {\rm st}\,X\;\llra\;\exists\,x\in V\,(\sta x\eqs X)$
for all $X\in\sV.$\vspace{1mm}

\noindent It is not assumed, in general, that $\eqs$ coincides
with the true equality on $\sV,$ but $\eqs$ has to be an
equivalence relation and satisfy the logic axioms for equality
with respect to $\ins$ and $\sta {\rm st}.$\vspace{5mm}

{\bf Ground model.} Thus it is assumed that there exist cardinals
$\vt$ such that $\bV\!_\vt$ is a model of ${\sf ZFC}.$ (This
assumption is valid during the proof of Theorem 2 below.) Let $\vt
$ be {\sf the least} among the cardinals of this kind. We use the
set $V=\bV\!_\vt$ as the ground ${\sf ZFC}$ model.

The first version of Theorem 1 is as follows:\vspace{3mm}

{\bf Theorem 2.} {\it Assume that the axiom of constructibility\/
$\bV=\bL$ holds. The model\/ $V$ has two regular\/
${\sf IST}$ extensions,\/ $\sV$ and\/ $\sV ',$ such that\/ $(\ast)$ is
false in\/ $\sV$ but true in\/ $\sV '.$}\vspace{3mm}

{\bf Proof.} The principal property implied by the minimality of
$\vt $ is that $\vt $ has countable cofinality. Let, indeed, $\vpi
_m (v_1,...,v_{n_m}),$ $m \in \omega ,$ be a recursive enumeration
of all parameterfree $\hspace{-1mm}\in\hspace{-1mm}$-formulas. It
is a theorem of ${\sf ZFC}$ (the Reflection principle, applied in
$V$) that for each integer $n$ there exists a cardinal $\kpa < \vt
$ such that $\bV\!_{\kpa }$ is an elementary submodel of $\bV\!_
\vt$
with respect to all sentences of type $ \vpi _m(x_1,...,x_{n_m}),$
where $m \leq n$ and $x \in \bV\!_\kpa .$ Let $\kpa _n$
denote the least cardinal $\kpa $ of such a kind; then $\kpa _n
\leq \kpa _{n+1}$ for all $n.$\vspace{2mm}

{\sf Lemma 3.} $\vt ={\rm sup}_{n\in\omega}\kpa _n.$\vspace{2mm}

{\sf Proof.} Let, on the contrary, $\vt >\kpa ={\rm
sup}_{n\in\omega}\kpa _n.$ By the definition of $\kpa _n,$ $V
'=\bV\!_\kpa $ is an elementary submodel of $V$ with respect
to {\it all\/} formulas $\vpi _m$ having sets in $V '$ as
parameters, hence a ${\sf ZFC}$ model, which contradicts the
choice of $\vt .$ $\Box$ \vspace{2mm}

The sequence of ordinals $\kpa _n$ plays an extremely important
role in the proof of Theorem 1 since it is the one that gives a
counterexample to $(\ast)$ in the extension of $V$ where $(\ast)$
fails.

The essential consequence of the assumption of $\bV=\bL$
here is that a certain relation $<_\bL$ wellorders the
universe of all sets $\bV$ in such a way that the following
property is guaranteed: given a cardinal $\vt $ such that $
\bV\!_\vt $ is a ${\sf ZFC}$ model, the relation $<_\bL$
wellorders $\bV\!_\vt$ with order type $\vt $ and is
$\hspace{-1mm}\in\hspace{-1mm}$-definable
in $\bV\!_\vt.$ This will be used in the construction of the model
$\sV ',$ where $(\ast)$ is true, and is irrelevant to the other
extension, $\sV.$\vspace{5mm}

{\sf The ultrafilter.} Both extensions, $\sV$ and $\sV ',$ are
constructed as ultrapowers of $V$ via a common ultrafilter, a kind
of {\it adequate} ultrafilters of Nelson [1977]. (Original
Nelson's construction includes infinite number of successive
ultrapowers; we show here that this can be managed an one-step
construction.) We introduce the {\it index set} \vspace{1mm}

\indent\indent $I = {\cal P}^{\rm fin}(V)=\{i\in V: i\hbox { is
finite}\}.$\vspace{1mm}

Let ${\rm Def}(V)$ denote the collection of all sets $X\subseteq
V,$ 1st order definable in $V$ by
$\hspace{-1mm}\in\hspace{-1mm}$-formulas having sets in $V$ as
parameters. \vspace{3mm}

{\sf Lemma 4.} {\it There exists an ultrafilter\/ $U$ over\/ $I$
satisfying the following two properties: \ \ {\rm (A)} \
$I_a=\{i \in I: a \in i\} \in U$ whenever $a \in V;$

{\rm (B)} \ $\{x \in V:$ {\rm the set} $P_x = \{i: \langle
i,x\rangle \in P\}$ {\rm is in} $U\}$ belongs to ${\rm
Def}(V)\hspace{-1mm}$\break
\indent\hspace{9mm} whenever $P \subseteq I\times V,$ $P \in
{\rm Def}(V).$ }\vspace{2mm}

{\sf Proof.} The construction is divided onto three
stages.

1. We define $U_0$ to be the collection of all sets of type $\{i
\in I: a \in i\},$ where $a \in V.$ It is evident that $U_0$ has
the finite intersection property (f.i.p.) which states that the
intersection of any finite number of sets contained in the
collection is nonempty.

2. We fix an enumeration $\chi_k(i,x),$ $k \geq 1,$ of all
parameterfree\break
$\in\hspace{-1mm}$-formulas with $i$
and $x$ as the only free
variables. Recall that $V$ is wellordered by the order relation
$<_\bL$ so that the order type of $V$ is $\vt .$ Let $x_\alpha
$ ($\alpha < \vt $) be the $\alpha $th element of $V$ with
respect to $<_\bL.$ The sequence $\langle x_\alpha : \alpha <
\vt \rangle $ belongs to ${\rm Def}(V)$ because $<_\bL$ restricted
to $V$ belongs to ${\rm Def}(V).$ We define \vspace{0mm}
$$
A_k(\alpha ) = \{i \in I: \chi_k(i,x_\alpha ) \;{\rm is\
true\ in\ } V\} \;\; {\rm and}\;\; C_k(\alpha ) = I \setminus
A_k(\alpha ).
$$
One can construct by induction on $k$ and, for a given $k,$ by
induction on $\alpha ,$ a collection of sets $T_k \subseteq
\vt ,$ $T_k \in {\rm Def}(V),$ such that the following sets
$$
\begin{array}{rll}
U_k\, = &\!\!\!\!\!\{A_k(\alpha ): \alpha \in T_k\}\!\!\!\!&\bigcup
\,\;\{C_k(\alpha ):\alpha\in\vt\setminus T_k\},\vspace{2mm}\\
U_{k\gamma } = &\!\!\!\!\!\{A_k(\alpha ): \alpha \in T_k \bigcap
\gamma \}\!\!\!\!&\bigcup\;\,\{C_k(\alpha ):\alpha\in\gamma\setminus
T_k\},\vspace{0mm}
\end{array}
$$
satisfy the condition that the union $U_{0} \bigcup ... \bigcup
U_{k-1}\bigcup U_{k\gamma}$ has the f.i. property for all $k \geq
1$ and $\gamma < \vt .$ The decision which of the sets
$A_k (\alpha),$ $C_k(\alpha)$ has to be adjoined to $U_k$ is
made so that we select $A_k(\alpha )$ provided this does not
violate f.i.p., and we select $C_k(\alpha )$ otherwise.

3. We set $U_{\infty} = \bigcup _{k\in\bN } U_k$ and extend
$U_{\infty}$ to an ultrafilter $U$ over $I.$ The ultrafilter $U$
is as required. One can easily verify (B) using the property
of definability of the sequence $\langle x_\alpha :\alpha
<\vt \rangle.$ $\Box $ \vspace{2mm}

It is assumed henceforth that $U$ is an ultrafilter given by the
lemma. Take notice that the property (B) of the ultrafilter $U$ is
essential only for the construction of the extension $\sV '$ but
not for $\sV.$

We introduce a convenient tool, the quantifier ``there exist
$U$-many" by \vspace{1mm}
$$
{\bf U}\,i\,\vpi (i) \;\;{\rm \ if\ and\ only\ if\
}\;\;\;\{i\in I:\vpi (i)\}\in U. \vspace{1mm}
$$
The following is the list of properties of ${\bf U}$ implied by
the definition of an ultrafilter and ( this regards (U5) and (U6))
the choice of the ultrafilter $U.$\vspace{2mm}

(U1) \ $\vpi\,\;\llra\,\; {\bf U}\, i\,\,\vpi
\;\;$ whenever $i$ is not free in $\vpi $;\vspace{1mm}

(U2) \ if $\forall\,i\;[\vpi (i)\;\lra\;\psi
(i)]$ then ${\bf U}\,i\;\vpi (i)\ \lra\ {\bf
U}\,i\;\psi (i)$;\vspace{1mm}

(U3) \ ${\bf U}\,i\;\vpi (i)\cj {\bf U}\,i\;\psi (i)\ \llra \ {\bf
U}\,i\;[\,\vpi (i)\cj\psi (i)\,]$;\vspace{1mm}

(U4) \ ${\bf U}\,i\;\neg \,\vpi (i)\ \llra\ \neg
\,{\bf U}\,i\;\vpi (i)$;\vspace{1mm}

(U5) \ if $a\in V$ then ${\bf U}\, i\,\, (a\in i)$;
\vspace{1mm}

(U6) \ Let $P\subseteq I\times V,$ $P\in {\rm Def} (V).$ Then
\vspace{-1mm}

\hfill $\{x\in V: {\bf U} \,i\;(\langle i,x\rangle \in P)\}\in {\rm
Def}(V).$\vspace{5mm}

{\sf The ``falsity" extension.} The union of ultrapowers of
$V$ via $U^r,$ $r\in\omega,$ is used to obtain a regular
${\sf IST}$ extension $\sV$ of $V$ where $(\ast)$ fails. We put
$$
\sV\hspace{-1mm}_r = \{f:f\ {\rm is\ a\ function},\ f: I^r\ \lra\
V\}.
$$
In particular, $\sV\hspace{-1mm}_0=\{\sta z:z\in V\},$ where $
\sta z = \{\langle 0,z\rangle\},$ since $I^0=\{\emptyset\}.$

The set $\sV = \bigcup_{r\in\omega }\sV\hspace{-1mm}_r$ is what we
call {\it the falsity extension}.

To continue notation, we let, for $F\in\sV,$ \ $r(F)$ denote the
unique\vspace{1mm} $r$ satisfying $F\in\sV\hspace{-1mm}_r.$ If
$F\in\sV,$ $q\geq r=r(F),$\vspace{1mm}
${\bf i} = \langle i_1,..., i_r, ..., i_q \rangle \in I^q,$
then we put $F[{\bf i}] = F(i_1, ..., i_r).$ Note that $F[{\bf i}]
= F({\bf i})$\vspace{1mm} whenever $r = q.$ We define
finally $\sta z[{\bf i}]=z$ for all $\sta z\in \sV\hspace{-1mm}_0$
and ${\bf i}\in I^{r}, \;r\geq 0.$

Let $F,G \in \sV$ and $r = {\rm max} \{r(F),r(G)\}.$ We set
\vspace{1mm}
$$
F\ins G\hspace{5mm}{\rm if\ and\ only\ if}\hspace{5mm}{\bf U}\,i_r
\;{\bf U}\,i_{r-1}... {\bf U}\,i_1 \;(F[{\bf i}] \in G[{\bf i}]);
$$
$$
F\eqs G\hspace{5mm}{\rm if\ and\ only\ if}\hspace{5mm}{\bf U}\,i_r
\;{\bf U} \,i_{r-1} ... {\bf U} \,i_1 \;(F[{\bf i}] = G[{\bf i}]);
\vspace{1mm}
$$
of course ${\bf i}$ denotes the sequence $i_1,...,i_r$.

The definition of standardness in $\sV$ is given by:
$$
\sta {\rm st}\,F\hspace{5mm} {\rm if\ and\ only\ if}\hspace{5mm}
{\rm there\ exists}\hspace{2mm} x \in V\hspace{2mm} {\rm such \
that}\hspace{2mm} F\eqs\sta x.
$$
So up to the relation $\eqs$ the level $\sV\hspace{-1mm}_0$ is just the
standard part of $\sV.$

Let, finally, $\Phi$ be a formula with parameters in $\sV.$
We define\break
$r(\Phi) = {\rm max}\{r(F):F$ occurs in $\Phi\}.$ If
in addition $r\geq r(\Phi)$ and ${\bf i}\in I^r,$ then let $\Phi
[{\bf i}]$ denote the result of replacing each $F$ that occurs in
$\Phi$ by $F[{\bf i}].$ Clearly $\Phi [{\bf i}]$ is a formula with
parameters in $V.$\vspace{4mm}

{\sf Proposition 5.} $\langle\sV; \eqs, \ins, \;\sta {\rm
st}\rangle$ {\it is a regular\/} ${\sf IST}$ {\it extension of\/}
$V$ {\it where\/ $(\ast)$ fails}.\vspace{2mm}

{\sf Proof.} The following principal statement plays the key
role.\vspace{3mm}

{\sf Lemma 6.} [\L o\' s Theorem] {\it Let\/ $\Phi$ be an internal
formula with parameters in\/ $\sV$ and suppose that\/ $r \geq r
(\Phi).$ Then }
$$
\Phi\ is\ true\ in\ \sV\;\;\llra\;\;{\bf U}\, i_r\,...\,{\bf
U}\,i_1\; (\Phi [i_1,..., i_r]\ \ is\ true\ in\ V).\vspace{0mm}
$$

{\sf Proof of the lemma}. The proof goes by induction on the
logical complexity of $\Phi.$ We abandon easy parts of the proof,
based on properties (U2), (U3), (U4) of the quantifier ${\bf U},$
and consider the induction step $\exists.$ Thus the lemma is to be
proved for a formula $\exists\,x\,\Phi (x)$ in the assumption that
the result holds for $\Phi (F)$ whenever $F \in\sV.$ We denote $r
= r (\Phi).$

The direction $\lra.$ Suppose that $\exists\,x\,\Phi (x)$ holds in
$\sV.$ Then $\Phi (F)$ holds in $\sV$ for some $F\in\sV.$ Let $p
= {\rm max}\,\{r,r(F)\}.$ To convert the reasoning into a more
convenient form, we let ${\bf i}$ and ${\bf j}$ denote sequences
$$
\langle i_1, ..., i_r\rangle\;\;(\in I^r)\hspace{3mm}{\rm and}
\hspace{3mm} \langle i_1, ..., i_r, ..., i_p\rangle\;\;(\in I^p)
$$
respectively. Further let ${\bf U\,i}$ and ${\bf U\,j}$ denote
sequences of quantifiers
$$
{\bf U}\,i_r\,...\,{\bf U}\,i_1\hspace{3mm}{\rm and}\hspace{3mm}
{\bf U}\,i_p\,...\,{\bf U}\,i_r\,...\,{\bf U}\,i_1.
$$
Thus ${\bf U\,j}\;\Phi (F)\,[{\bf j}]$ holds by the induction
hypothesis. We note that, for all\vspace{1mm}\break
{\bf j}, $\Phi (F)\,[{\bf j}]\ \lra\ \exists\,x\,\Phi (x)\,[{\bf
j}].$ Hence ${\bf U \,j}\;\exists\,x\,\Phi (x)\,[{\bf j}]$ is true
by (U2). We note\vspace{1mm}\break
also that the formula $\exists\,x\,\Phi (x)\,[{\bf j}]$ coincides
(graphically) with $\exists\,x\,\Phi (x)\,[{\bf
i}]$\vspace{1mm}\break
because $r(\exists\,x\,\Phi (x)) = r \leq p.$ Hence, deleting the
superfluous quantifiers\vspace{1mm}\break
by (U1), we obtain ${\bf U\,i}\;\exists\,x\; \Phi (x)\,[{\bf
i}].$

The direction $\longleftarrow.$ Let $\Phi (x)$ be $\Phi
(x,G,H,...,),$ where $G,H,... \in\sV.$ Suppose that ${\bf
U\,i}\;\exists\,x\;\Phi (x)\,[{\bf i}]$ holds, that is,
$$
{\bf U\,i}\;[\,\exists\,x\;\Phi (x,G[{\bf i}],H[{\bf i}],...)
\hspace{2mm} {\rm is\ true\ in\ } V\,].
$$
For each ${\bf i}\in I^r,$ if there exists some $x \in V$ such
that $\Phi (x,G[{\bf i}],H[{\bf i}],...)$ is true in $V,$
then we let $F({\bf i})$ be one of $x$ of such kind; otherwise let
$F({\bf i}) = \emptyset.$ By definition, $F\in \sV\hspace{-1mm}_r$ and
$$
\forall\,{\bf i}\in I^r\;[\,\exists\,x\;\Phi (x)\,[{\bf
i}]\ \lra\ \Phi (F)\,[{\bf i}]\,],
$$
therefore ${\bf U\,i}\;\exists\,x\;\Phi (x)\,[{\bf i}]\ \lra\ {\bf
U\,i}\;\Phi (F)\,[{\bf i}]$ by (U2). Recall that the left-hand
side of the last implication has been supposed to be true. So the
right-hand side is also true. Then $\Phi (F)$ holds in $\sV$ by
the induction hypothesis, and we are done. $\Box$\vspace{2mm}

The just proved lemma easily implies logical equality axioms for
$=,$ and Transfer, therefore all ${\sf ZFC},$ in $\sV.$
Standardization is evident because every set $V$ of the form
$V=\bV\!_\vt$ has the property that if $Y\subseteq
X\in V$ then $Y\in V.$ We prove Idealization. \vspace{1mm}

Thus let $\vpi (x,a)$ be an internal formula with parameters
in $\sV.$ We denote $r = r (\vpi)$ and prove the following:
$$
\fstf A\;\exists\,x\;\forall\,a \in A\;\vpi (x,a)\ \lra\
\exists\,x\;\fst a\;\vpi (x,a)
$$
in $\sV.$ (The implication $\longleftarrow$ does not need a
special consideration because it follows from Standardization that
elements of finite standard sets are standard, see Nelson [1977].)
Lemma 6 converts the left-hand side to the form:
$$
\forall\hspace{0.5pt}^{\rm fin}A\;\subseteq\,V\;{\bf
U}\,i_r\,...\,{\bf U}\,i_1\; \exists\,x\;\forall \,a
\in A\; (\vpi (x,a)\,[i_1, ..., i_r]).
$$
Recall that $I$ consists of all finite subsets of $V,$ so we may
replace the\vspace{1mm} variable
$A$ by $i,$ having in mind that $i \in I.$ Further\vspace{1mm}
define $\tilde A: I^{r+1}\ \lra\ V$ by
$\tilde A(i_1,..., i_r, i) = i.$ Then $\tilde A \in\sV.$ The
left-hand
side takes the form
$$
\forall\,i\;{\bf U}\, i_r\,...\,{\bf U}\, i_1 \; (\exists\,x\;
\forall\,a\in \tilde A\; \vpi (x,a))\,[i_1,..., i_r, i].
$$
Changing $\forall\,i$ by ${\bf U}\,i,$ we obtain $\exists\,x\;
\forall\, a \in \tilde A\; \vpi (x,a)$ in $\sV$ again by the
lemma. So, to verify the right-hand side of Idealization, it
suffices to prove $\sta a \in \tilde A$ in $\sV$ for all $a \in
V.$ This is equal to
$$
{\bf U}\,i\;{\bf U}\,i_r\,...\,{\bf U}\,i_1\;(a \in \tilde A
\,[i_1, ... ,i_r, i]),
$$
by the lemma, and then to ${\bf U} \,i\;{\bf U} \,i_r
\,...\,{\bf U}\, i_1\; (a \in i)$ by the definition of $\tilde A.$
So apply (U1) and complete the proof of Idealization in $\sV.$
\vspace{1mm}

Thus $\sV$ is an ${\sf IST}$ model. Moreover it is a regular
extension of $V$: one can easily verify the required properties of
the embedding $^\ast.$ To complete the proof of Proposition 5 it
remains to show that $(\ast)$ does not hold in $\sV.$\vspace{1mm}

We use the sequence of ordinals $\kpa _n.$ Let $F\in \sV\!_0$ be
defined by
$$
F(i)=\{\langle n,\kpa _n\rangle : \langle n,\kpa _n\rangle \in
i\} \;\;\;{\rm for\ all\ } i\in I.
$$
It is true in $\sV$ by Lemma 6 that $F$ is a function defined on a
subset of integers, and, for every $n\in \omega ,$ it is also true
in $\sV$ that $F(^{\ast}\!n)$ is defined and equal to $\sta\kpa
_n,$ hence standard. Thus the left-hand side of $(\ast)$ is
satisfied by $F.$

The right-hand side cannot be satisfied since it would imply that
there exists $g\in V$ s.t. $g(n)=\kpa _n$ for all $n,$ a
contradiction with Lemma 3. $\Box$\vspace{5mm}

{\sf The ``truth" extension.} We continue the proof of Theorem 2.
To get rid of such elements of $\sV$ as the considered above $F,$
we build up the required ultrapower using only {\it definable\/}
functions. Thus we set
$$
\sV '\hspace{-1mm}_r=\{F:I^r\;\lra\;V, \;F\in {\rm
Def}(V)\}\hspace{3mm}{\rm for\ all\ }r.
$$
The model $\sV '=\bigcup_{r\in \omega }\,\sV '\hspace{-1mm}_r$ is
the {\it truth extension}.

All other relevant definitions are the same as above.\vspace{3mm}

{\sf Proposition 7.} {\it $\sV '$ is a regular\/ ${\sf IST}$
extension of\/ $V$ where $(\ast)$ is true.}\vspace{2mm}

{\sf Proof.} To prove that $\sV '$ is an ${\sf IST}$ model and a
regular extension of $V$ one can copy the proof of Proposition 5.
The only notable difference is related to the direction
$\longleftarrow$ in the proof of Lemma 6. The problem is that
$F\in {\rm Def}(V)$ should be guaranteed; otherwise one
cannot claim that $F\in\,\sV '.$

To fix the problem we define $F({\bf i})$ to be the $<_\bL
\hspace{-1mm}$-least $x$ satisfying the corresponding property.
Thus $F\in {\rm Def}(V)$ because $<_\bL$ is definable in
$V.$\vspace{1mm}

It remains to prove that $(\ast)$ is true in $\sV '.$ Thus let
$F\in\sV '$ be such that it is true in $\sV '$ that $F$ is a
function defined on integers and $F(n)$ is standard for every
standard $n.$ We set, for all $n\in \omega $ and $x\in V,$
$$
g(n)=x\;\;\;{\rm if\ and\ only\ if\ }\;\;\;F(\sta
n)=\sta x\;\;{\rm in}\;\sV '.
$$
Then $g$ is defined at all $n\in \omega $ and $F(\sta n)=\sta
g(n)$ in $\sV '$ for all $n.$

Moreover $g$ is definable in $V.$ Indeed,
$$
g(n)=x\,\;\llra\;\,{\bf U}\,i_r\,...\,{\bf U}\,i_1\;
[\,(F(\sta n)=\sta x)\,[i_1,...,i_r]\,]
$$
by Lemma 6 ($r=r(F)$). The relation in outer brackets (in $r+2$
variables) is definable in $V$ since $F\in {\rm Def}(V).$ Finally
the action of ${\bf U}$ keeps the definability by (U6). Thus $g$
is definable in $V,$ therefore $g\in V.$

To end the proof we define $G=\sta g.$ Then, for all $n,$ it is
true in $\sV '$ that $F(\sta n)=G(\sta n).$ $\Box$ \vspace{2mm}

This completes the proof of Theorem 2. $\Box$ \vspace{5mm}

{\bf The proof without models.} To avoid the assumption of the
existence of a cardinal $\vt $ such that $\bV\!_\vt$ is a
${\sf ZFC}$ model, we apply a logical trick. We extend the
$\hspace{-1mm}\in \hspace{-1mm}$-language of ${\sf ZFC}$ by a
special constant $\dot V $ and add
the axiom of constructibility $\bV=\bL,$ the statement
\vspace{1mm}

({\rm Mod}) \ \ $\dot V$ has the form $\dot V=\bV\!_\vt $
for a cardinal $\vt $ of countable cofinality;\vspace{1mm}

\noindent and the collection of all sentences of type: ``$A$ holds
in $\dot V$", where $A$ is an axiom of ${\sf ZFC},$ to the list of
${\sf ZFC}$ axioms. Let us denote the extension by ${\sf ZFC}\,[\dot
V].$ Thus
$$
{\sf ZFC}[\dot V] = {\sf ZFC} + ({\rm Mod}) + [\bV=\bL]
+\{A^{\dot V}\!:A {\rm \ is\ a\ {\sf ZFC}\ axiom}\},
$$
where $A^{\dot V}$ is the relativization of $A$ to $\dot V.$ (To
obtain $A^{\dot V}$ one has to replace every quantifier ${\rm
Q}\,x$ in $A$ by ${\rm Q}\,x\in \dot V.$)\vspace{3mm}

{\sf Proposition 8.} {\it ${\sf ZFC}$ and\/ ${\sf ZFC}\,[\dot V]$
are equiconsistent.}\vspace{2mm}

{\sf Proof.} It is sufficient to prove that an arbitrary {\it
finite\/} extension of ${\sf ZFC}$ of the type
$$
{\sf ZFC}^{\rm fin}[\dot V] = {\sf ZFC} + ({\rm Mod}) + [
\bV=\bL] +\{A_1^{\dot V},...,A_n^{\dot V}\},
$$
where $A_1,...,A_n$ are ${\sf ZFC}$ axioms, is equiconsistent
with ${\sf ZFC} + [\bV=\bL].$

We use the following statement (a kind of Reflection principle).
For any finite list $A_1,...,A_n$ of
$\hspace{-1mm}\in\hspace{-1mm}$-sentences it is a
theorem of ${\sf ZFC}$ that there exists a cardinal $\vt $ of
countable cofinality such that
$\bV\!_\vt$ is an elementary submodel of the universe
$\bV$ of all sets with respect to $A_1,...,A_n.$ In
particular, since all of $A_1,...,A_n$ are ${\sf ZFC}$ axioms
here, all of them are true in $\bV\!_\vt.$

We define, arguing in ${\sf ZFC}+[\bV=\bL],$ $\vt $
to be the least among such ordinals and obtain an interpretation
of ${\sf ZFC}^{\rm fin}[\dot V]$ in ${\sf ZFC}+[\bV=\bL]$ by
interpreting $\dot V$ as $\bV\!_\vt .$ $\Box$\vspace{3mm}

{\bf Proof of Theorem 1.} Let, on the contrary, $\Phi $ be a
(parameterfree)\break
$\in\hspace{-1mm}$-sentence such that the equivalence
$\bS\;\llra\;(\ast)$ is an ${\sf IST}$ theorem. Let ${\sf IST}^{\rm
fin}$ be a finite part of ${\sf IST}$ sufficient to prove the
equivalence.

We argue in ${\sf ZFC}\,[\dot V].$ By axiom ({\rm Mod}), $\dot V$
has the form $\dot V=\bV\!_\vt ,$ where $\vt $ is a
cardinal of countable cofinality. Let $\langle \kpa _n:n\in
\omega \rangle $ be\vspace{1mm} a cofinal in $\vt $
countable sequence of ordinals. Let, finally, $\sta {\dot V}$ and
$\sta\dot V '$ be the ``falsity" and ``truth" extensions of $\dot
V$ defined as above.

It is not assumed, of course, that $\dot V$ is a model of ${\sf
ZFC}.$ However the given above reasoning related to ``falsity" and
``truth" extensions can be converted to the form: given an axiom
$A$ of ${\sf IST},$ there exists a finite list $Z_1,...,Z_n$ of
${\sf ZFC}$ axioms such that $A$ is true in $\sV$ and $\sV '$
provided $V=\bV\!_\vt $ is a model of $Z_1,...,Z_n.$ Therefore
both $\sta\dot V$ and $\sta\dot V '$ are models of ${\sf IST}^{\rm
fin}.$

By the same argument, $(\ast)$ is false in $\sta\dot V$ and true
in $\sta\dot V '.$ Thus $\Phi $ is false in $\sta\dot V$ and true
in $\sta\dot V '$ by the choice of ${\sf IST}^{\rm fin}.$

By the same argument, this implies that $\Phi $ is simultaneously
false and true in $\dot V.$ This proves that ${\sf ZFC}\,[\dot V]$
is inconsistent, a contradiction with Lemma 8 and the assumption
that {\sf ZFC} is consistent. $\Box$ \vspace{4mm}

\begin{center}
{\bf References}
\end{center}

\noindent I. van den Berg [1987]\\
\indent {\it Nonstandard asymptotic analysis} (Lecture Notes in
Math. 1249, Springer).\vspace{1mm}

\noindent F.Diener and M.Diener [1988]\\
\indent Some asymptotic results in ordinary differential
equations. in: N.Cutland\break\indent
(ed.) {\it Nonstandard analysis and its applications} \
(London Math. Soc.\break\indent
Student Texts 10, Cambridge Univ. Press), pp. 282 --
297.\vspace{1mm}

\noindent F.Diener and G.Reeb [1989]\\
\indent {\it Analyse non standard\/} (Herrmann
Editeurs).\vspace{1mm}

\noindent V.G.Kanovei [1991]\\
\indent Undecidable hypotheses in Edward Nelson's internal set
theory, {\it Russian\break\indent
Math. Surveys}, 46, pp. 1 -- 54.\vspace{1mm}

\noindent E.Nelson [1977]\\
\indent Internal set theory; a new approach to nonstandard
analysis, {\it Bull. Amer.\break\indent
Math. Soc.} 83, 1165 -- 1198.\vspace{1mm}

\noindent E.Nelson [1988]\\
\indent The syntax of nonstandard analysis, {\it Ann. Pure and
Appl. Log.} 38, 123\break\indent -- 134.

\noindent M.Reeken [1992]\\
\indent On external constructions in internal set theory. {\it
Expositiones Mathe-\break\indent
maticae,\/} 10, 193 -- 247.

\end{document}